\documentclass[reqno]{amsart}

\usepackage{amsmath}
\usepackage{amsthm}
\usepackage{amsfonts}
\usepackage{amssymb}
\usepackage{latexsym}
\usepackage{amscd}
\usepackage{stmaryrd}
\usepackage{cite}
\usepackage{color}

\theoremstyle{plain}
\newtheorem{thm}{Theorem}
\newtheorem{cor}[thm]{Corollary}

\newfont{\scyr}{wncyr10 scaled 550}

\definecolor{darkgreen}{rgb}{0.0625,0.64,0.0625}
\usepackage[%dvipdfm,
            pdfstartview=FitH,
            CJKbookmarks=true,
            bookmarksnumbered=true,
            bookmarksopen=true,
            colorlinks,
            pdfborder=001,
            linkcolor=blue,
            anchorcolor=green,
            citecolor=red
            ]{hyperref}

\allowdisplaybreaks

\begin{document}

\title{Derivation relations and duality for the sum of multiple zeta values}

\date{\today\thanks{The author is supported by the National Natural Science Foundation of
China (Grant No. 11471245) and the Natural Science Foundation of Shanghai (Grant No. 14ZR1443500).} }

\author{Zhonghua Li}

\address{School of Mathematical Sciences, Tongji University, No. 1239 Siping Road,
Shanghai 200092, China}

\email{zhonghua\_li@tongji.edu.cn}



\begin{abstract}
We show that the duality relation for the sum of multiple zeta values with fixed weight, depth and $k_1$ is deduced from the derivation relations, which was first conjectured by N. Kawasaki and T. Tanaka.
\end{abstract}

\maketitle

%%---------------------------------------------------------------------------
%%------------------------Content-------------------------------------------
%%----------------------------------------------------------------------------

For positive integers $n,k_1,k_2,\ldots,k_n$ with $k_1\geqslant 2$, the multiple zeta value is defined by the convergent series
$$\zeta(k_1,k_2,\ldots,k_n)=\sum\limits_{m_1>m_2>\cdots>m_n>0}\frac{1}{m_1^{k_1}m_2^{k_2}\cdots m_n^{k_n}}.$$
We call $k_1+\cdots+k_n$ the weight and $n$ the depth, respectively.

In the famous paper \cite{Ihara-Kaneko-Zagier}, K. Ihara, M. Kaneko and D. Zagier gave the framework of the extended double shuffle relations and conjectured that the extended double shuffle relations give all linear relations among multiple zeta values. So it is interesting to see whether all other linear relations of multiple zeta values can be deduced from the extended double shuffle relations or not. For some relative work on this direction, one may refer to \cite{Ihara-Kaneko-Zagier,Kajikawa,Kawasaki-Tanaka,Li2010,Li2013,Li-Qin}.

This short note is inspired by the recent paper \cite{Kawasaki-Tanaka} of N. Kawasaki and T. Tanaka. They showed that the duality relations of the following two cases are implied by the derivation relations, which are contained in the extended double shuffle relations:
\begin{itemize}
  \item the case for the double zeta values,
  \item the case for the sum of multiple zeta values with fixed weight, depth and $k_1=2$.
\end{itemize}
More generally, they conjectured that the duality relation for the sum with fixed weight, depth and $k_1$ is subject to the derivation relations. We prove this conjecture here.

To state the result, we use the algebraic setup introduced by M. E. Hoffman \cite{Hoffman1997}. Let $\mathfrak{h}=\mathbb{Q}\langle x,y\rangle$ be the noncommutative polynomial algebra generated by $x$ and $y$ over the filed $\mathbb{Q}$ of rational numbers. Let $\mathfrak{h}^0=\mathbb{Q}+x\mathfrak{h}y$ be a subalgebra. We define the $\mathbb{Q}$-linear map $Z:\mathfrak{h}^0\rightarrow \mathbb{R}$ by $Z(1)=1$ and
$$Z(x^{k_1-1}y\cdots x^{k_n-1}y)=\zeta(k_1,\ldots,k_n),$$
where $\mathbb{R}$ is the field of real numbers, $n,k_1,\ldots,k_n$ are positive integers with $k_1\geqslant 2$.

Now we state the duality and the derivation relations of multiple zeta values. Let $\tau$ be the anti-automorphism of the algebra $\mathfrak{h}$ determined by
$$\tau(x)=y,\qquad \tau(y)=x.$$
Then the duality relations \cite{Zagier} are nothing but
$$(1-\tau)(\alpha)\in\ker Z$$
for any $\alpha\in\mathfrak{h}^0$. A derivation of $\mathfrak{h}$ is a $\mathbb{Q}$-linear map $D:\mathfrak{h}\rightarrow \mathfrak{h}$ satisfying the Leibniz rule
$$D(\alpha\beta)=D(\alpha)\beta+\alpha D(\beta)\qquad (\forall \alpha,\beta\in\mathfrak{h}).$$
For a positive integer $n$, let $\partial_n$ be the derivation of $\mathfrak{h}$ determined by
$$\partial_n(x)=x(x+y)^{n-1}y,\qquad \partial_n(y)=-x(x+y)^{n-1}y.$$
Then the derivation relations \cite{Ihara-Kaneko-Zagier} are just
$$\partial_n(\alpha)\in\ker Z$$
for any positive integer $n$ and any $\alpha\in\mathfrak{h}^0$. For a formal variable $u$, let
$$\Delta_u=\exp\left(\sum\limits_{n=1}^\infty\frac{\partial_n}{n}u^n\right),$$
which is an automorphism of $\mathfrak{h}[[u]]$ and satisfies
$$\Delta_u(x)=x\frac{1}{1-yu},\quad \Delta_u(y)=(1-xu-yu)\frac{y}{1-yu},\quad \Delta_u(x+y)=x+y.$$
As in \cite{Ihara-Kaneko-Zagier},  the derivation relations can be restated as
$$(\Delta_u-1)(\alpha)\in\ker Z$$
for any $\alpha\in\mathfrak{h}^0$.

Using the algebraic notations, the conjecture of N. Kawasaki and T. Tanaka in \cite{Kawasaki-Tanaka} can be stated as
\begin{align}
(1-\tau)\left(\sum\limits_{a_1+\cdots+a_{l-1}=k-m-l\atop a_1,\ldots,a_{l-1}\geqslant 0}x^myx^{a_1}y\cdots x^{a_{l-1}}y\right)\overset{?}{\in}\sum\limits_{n\geqslant 1}\partial_n(\mathfrak{h}^0)
\label{Eq:Duality}
\end{align}
for any positive integers $k,m,l$ with $k\geqslant m+l$. We compute the generating function of the left-hand side of the above conjectured relation \eqref{Eq:Duality}. Let $u,v,w$ be independent formal variables, which commute with each other as well as with $x$ and $y$. We have
\begin{align*}
&\sum\limits_{k,m,l\geqslant 1\atop k\geqslant m+l}\left(\sum\limits_{a_1+\cdots+a_{l-1}=k-m-l\atop a_1,\ldots,a_{l-1}\geqslant 0}x^myx^{a_1}y\cdots x^{a_{l-1}}y\right)u^{m-1}v^{l-1}w^{k-m-l}\\
=&\frac{x}{1-xu}y\frac{1}{1-xw-yv}(1-xw)\\
=&\frac{x}{1-xu}y+\frac{x}{1-xu}y\frac{1}{1-xw-yv}yv\in\mathfrak{h}^0[[u,v,w]].
\end{align*}
Then the generating function of the left-hand side of \eqref{Eq:Duality} is
$$\frac{x}{1-xu}y-x\frac{y}{1-yu}+\left(\frac{x}{1-xu}y\frac{1}{1-xw-yv}y-x\frac{1}{1-xv-yw}x\frac{y}{1-yu}\right)v,$$
which is represented by the maps $\Delta_u,\Delta_v$ and $\Delta_w$ as in the following main theorem.

\begin{thm}\label{Thm:MainThm}
In $\mathfrak{h}^0[[u,v,w]]$, we have
\begin{align}
\frac{x}{1-xu}y-x\frac{y}{1-yu}=(1-\Delta_u)\left(\frac{x}{1-xu}y\right)
\label{Eq:Duality-Zeta}
\end{align}
and
\begin{align}
&\frac{x}{1-xu}y\frac{1}{1-xw-yv}y-x\frac{1}{1-xv-yw}x\frac{y}{1-yu}\nonumber\\
=&\frac{1}{v-w}(\Delta_v-\Delta_w)\left(x\frac{1}{1-xu-xv+(x^2+yx)uv}y\frac{1}{1-xw}(1-xw-yw)\right)\nonumber\\
&+(1-\Delta_u)\left(x\frac{1}{1-xu-xv+(x^2+yx)uv-yw}(1-xu-yu)\frac{x}{1-xu}y\right).
\label{Eq:Duality-k1}
\end{align}
\end{thm}

Note that taking $u=0$ and $v=0$ respectively in \eqref{Eq:Duality-k1}, we get \cite[Theorem]{Kawasaki-Tanaka} in terms of generating functions.
And by comparing the coefficients of the monomial $u^{m-1}v^{l-1}w^{k-m-l}$, we prove the conjecture of N. Kawasaki and T. Tanaka in \cite{Kawasaki-Tanaka}.

\begin{cor}
For any positive integers $k,m,l$ with $k\geqslant m+l$, we have
\begin{align*}
(1-\tau)\left(\sum\limits_{a_1+\cdots+a_{l-1}=k-m-l\atop a_1,\ldots,a_{l-1}\geqslant 0}x^myx^{a_1}y\cdots x^{a_{l-1}}y\right)\in\sum\limits_{n\geqslant 1}\partial_n(\mathfrak{h}^0).
\end{align*}
\end{cor}

The above corollary shows that the duality relation for the sum of multiple zeta values with fixed weight, depth and $k_1$ can be deduced from the derivation relations.
For example, taking $k=m+l$, we find the duality
$$\zeta(m+1,\underbrace{1,\ldots,1}_{l-1})=\zeta(l+1,\underbrace{1,\ldots,1}_{m-1})$$
is implied from the derivation relations.

Finally, we give a proof of Theorem \ref{Thm:MainThm}.

\noindent {\bf Proof of Theorem \ref{Thm:MainThm}.}
Since
$$\Delta_u(1-xu)=1-xu\frac{1}{1-yu}=(1-xu-yu)\frac{1}{1-yu},$$
we get \eqref{Eq:Duality-Zeta}.

Since
$$\Delta_u(1-xu-xv+(x^2+yx)uv-yw)=(1-xu-yu)(1-xv-yw)\frac{1}{1-yu},$$
we find
\begin{align*}
&\Delta_u\left(x\frac{1}{1-xu-xv+(x^2+yx)uv-yw}(1-xu-yu)\frac{x}{1-xu}y\right)\\
=&x\frac{1}{1-xv-yw}x\frac{y}{1-yu}.
\end{align*}
Similarly, we have
$$\Delta_v(1-xu-xv+(x^2+yx)uv)=(1-xv-yv)(1-xu)\frac{1}{1-yv},$$
which implies that
\begin{align*}
&\Delta_v\left(x\frac{1}{1-xu-xv+(x^2+yx)uv}y\frac{1}{1-xw}(1-xw-yw)\right)\\
=&\frac{x}{1-xu}y\frac{1}{1-xw-yv}(1-xw-yw)\\
=&\frac{x}{1-xu}y+(v-w)\frac{x}{1-xu}y\frac{1}{1-xw-yv}y.
\end{align*}
And since
$$\Delta_w(1-xu-xv+(x^2+yx)uv)=(1-xu-xv+(x^2+yx)uv-yw)\frac{1}{1-yw},$$
we have
\begin{align*}
&\Delta_w\left(x\frac{1}{1-xu-xv+(x^2+yx)uv}y\frac{1}{1-xw}(1-xw-yw)\right)\\
=&x\frac{1}{1-xu-xv+(x^2+yx)uv-yw}(1-xw-yw)y.
\end{align*}
Hence the right-hand side of \eqref{Eq:Duality-k1}  is
\begin{align*}
&\frac{1}{v-w}\frac{x}{1-xu}y+\frac{x}{1-xu}y\frac{1}{1-xw-yv}y\\
&-\frac{1}{v-w}x\frac{1}{1-xu-xv+(x^2+yx)uv-yw}(1-xw-yw)y\\
&+x\frac{1}{1-xu-xv+(x^2+yx)uv-yw}(1-xu-yu)\frac{x}{1-xu}y\\
&-x\frac{1}{1-xv-yw}x\frac{y}{1-yu}.
\end{align*}
Direct computation shows that
\begin{align*}
&(v-w)(1-xu-yu)\frac{x}{1-xu}-(1-xw-yw)\\
=&\left[(v-w)(1-xu-yu)x-(1-xw-yw)(1-xu)\right]\frac{1}{1-xu}\\
=&-(1-xu-xv+(x^2+yx)uv-yw)\frac{1}{1-xu}.
\end{align*}
Then it is easy to see that the right-hand side of \eqref{Eq:Duality-k1} is just the left-hand side of \eqref{Eq:Duality-k1}.
\qed

%%-----------------------------------------------------------------------
%%----------------------------References----------------------------------
%%-----------------------------------------------------------------------

\end{document}